\documentclass[11pt]{amsart}

\usepackage[usenames,dvipsnames]{xcolor}
\usepackage{amsmath,amsthm,amssymb,amsfonts,amscd, graphicx, enumerate, enumitem, hyperref, dsfont, bbm}
\usepackage{stmaryrd}
\usepackage{graphics}
\usepackage{tikz-cd}
\usepackage{stackrel}
\usepackage{mathtools}
\usepackage{wrapfig}
\usepackage[cal=boondoxo]{mathalfa}
\usepackage[font=footnotesize]{caption}
\usepackage{extarrows,multirow,colortbl,hhline}
\usepackage{array,diagbox}
\usepackage{parskip}
\usepackage[numbers]{natbib}
\setcitestyle{open={[},close={]}}
\usepackage{import}
\usepackage{pdfpages}
\usepackage{transparent}

\usepackage{slashed}
\usepackage{marginnote}
\usepackage{centernot}

\makeatletter 
\newsavebox{\@brx}
\newcommand{\llangle}[1][]{\savebox{\@brx}{\(\m@th{#1\langle}\)}%
  \mathopen{\copy\@brx\mkern2mu\kern-0.9\wd\@brx\usebox{\@brx}}}
\newcommand{\rrangle}[1][]{\savebox{\@brx}{\(\m@th{#1\rangle}\)}%
  \mathclose{\copy\@brx\mkern2mu\kern-0.9\wd\@brx\usebox{\@brx}}}
\makeatother

\hypersetup{colorlinks=true, linkcolor=MidnightBlue, citecolor=MidnightBlue, urlcolor=black}

\definecolor{IBMyellow}{RGB}{255,176,0}
\definecolor{IBMorange}{RGB}{254,97,0}
\definecolor{IBMred}{RGB}{220,38,127}
\definecolor{IBMpurple}{RGB}{120,94,240}
\definecolor{IBMblue}{RGB}{100,143,255} 

\usepackage[margin=1.3in,headsep=.3in, top=1.2in, bottom=1.15in]{geometry}

\newtheoremstyle{plain}
  {5pt plus3pt minus3pt}   
  {5pt plus3.3pt minus3.3pt}   
  {\slshape}  
  {0pt}       
  {\bfseries} 
  {.}         
  {5pt plus 1pt minus 1pt} 
  {}          

\newtheorem{theorem}{Theorem}[section]
\newtheorem{corollary}[theorem]{Corollary}
\newtheorem{lemma}[theorem]{Lemma}
\newtheorem{proposition}[theorem]{Proposition}
 
\newtheorem{question}[theorem]{Question}
\newtheorem*{question*}{Question} 
 
\newtheorem{definition}[theorem]{Definition}
\newtheorem{remark}[theorem]{Remark}

\newtheoremstyle{example}
  {2pt plus3.3pt minus3.6pt}   
  {2pt plus3.3pt minus3.3pt}   
  {}  
  {0pt}       
  {\bfseries} 
  {.}         
  {5pt plus 1pt minus 1pt} 
  {}          
\theoremstyle{example}

\theoremstyle{definition}

\DeclareMathOperator{\Int}{int}
\DeclareMathOperator{\rank}{rank}

\newcommand{\CC}{\mathbb{C}}
\newcommand{\RR}{\mathbb{R}}
\newcommand{\ZZ}{\mathbb{Z}}

\newcommand{\PP}{\mathbb{P}}
\newcommand{\id}{\operatorname{id}}

\newcommand{\homeo}{\underset{\text{top}}{\cong}}
\newcommand{\diffeo}{\underset{\text{sm}}{\cong}}

\newcommand{\MCG}{\mathrm{MCG}}
\newcommand{\R}{\mathcal{R}}
\newcommand{\FBF}{\mathrm{FBF}}
\newcommand{\FSW}{\mathrm{FSW}}

\newcommand{\s}{\mathfrak{s}}
\newcommand{\SW}{\mathrm{SW}}
\newcommand{\sw}{\mathrm{sw}}
\newcommand{\G}{\mathcal{G}}
\newcommand{\BF}{\mathrm{BF}}

\newcommand{\HH}{\mathbb{H}}

\newcommand{\Sp}{\operatorname{Sp}}

\newcommand{\U}{\operatorname{U}}

\newcommand{\Pin}{\mathrm{Pin}}
\newcommand{\DD}{\mathbb{D}}
\newcommand{\std}{\mathrm{std}}

\usepackage{fancyhdr}
\usepackage[affil-it]{authblk} 
\usepackage{etoolbox}
\usepackage{lmodern}
\makeatletter
\patchcmd{\@maketitle}{\LARGE \@title}{\fontsize{16}{19.2}\selectfont\@title}{}{}
\makeatother

\renewenvironment{proof}{{\sc Proof:}}{\hfill\ensuremath{\blacksquare}}

\author[A.\ Shivkumar]{Abhishek Shivkumar}
\address{Department of Mathematics, University of Texas, Austin, TX 78712}
\email{ashivkum@utexas.edu}

\title[Exotic $\RR^4$'s with Compactly Supported Diffeomorphisms]{\large Exotic $\RR^4$'s with Compactly Supported Diffeomorphisms}

\begin{document}
\vspace*{-.375in}

\bigskip

\begingroup
\def\uppercasenonmath#1{}
\let\MakeUppercase\relax
\maketitle
\endgroup
\thispagestyle{empty}

\vspace{-.25in}

\begin{center}\small
\textsc{Abhishek Shivkumar}
\end{center}


\bigskip

\begin{center} \begin{minipage}{.88\linewidth}\footnotesize

\textsc{Abstract.} We show that for diffeomorphisms of closed simply-connected $4$-manifolds dissolving under a single stabilization, there is a natural induced diffeomorphism on an exotic copy of $\RR^4$. We thereby construct uncountably many exotic copies of $\RR^4$ with nontrivial compactly supported mapping class groups of arbitrarily large rank. This follows from a modification of the construction of the diffeomorphism corks of \cite{krushkal2024corks} that makes their exteriors simply-connected. 

\end{minipage}
\end{center}

\bigskip

\section{Introduction}\label{sec:intro}

\smallskip

A diffeomorphism of a $4$-manifold is \emph{exotic} if it is topologically but not smoothly isotopic to the identity. Parallel to the story for exotic smooth structures on (closed, smooth, simply-connected) $4$-manifolds, the topological classification of automorphisms of such manifolds is governed by algebraic topology, and gauge theory may be used to obstruct isotopy to the identity. Exoticness in both settings dissolves after a finite number of stabilizations by $S^2 \times S^2$. 

In analogy to the classical cork theorem for exotic pairs, it has recently been shown that exotic diffeomorphisms $f:X \rightarrow X$ which are smoothly isotopic to the identity after a single stabilization may be localized to a contractible submanifold $C$. By localize, here, we mean that $f$ may be isotoped to $f'$ which is the identity on $X \setminus C$, so $f'$ is supported entirely on $C$. Such a manifold $C$ is called a \emph{diffeomorphism cork} for $f$ \cite[Theorem 4.1]{krushkal2024corks} (quoted below as Theorem~\ref{thm:4.1}). In fact, we can produce diffeomorphism corks for any finite collection of diffeomorphisms as above \cite[Theorem 4.2]{krushkal2024corks} (quoted below as Theorem~\ref{thm:4.2}).

Classical corks can be constructed to have simply-connected exterior. In this paper, we extend the above results by proving an analogous result for diffeomorphism corks: 
\begin{proposition}\label{prop:intro}
  Let $\{f_i\}_{i=1}^n$ be any finite collection of self-diffeomorphisms of a closed, simply-connected, smooth $4$-manifold $X$, each $1$-stably isotopic to the identity. There exists a smooth, compact, contractible submanifold $C \subset X$ and a collection of self-diffeomorphisms $\{f_i'\}_{i=1}^n$ such that $f_i'$ is smoothly isotopic to $f_i$, and the $f_i'$ are all supported on $C$. Moreover, $C$ satisfies $\pi_1(X\setminus C) = 1$. 
\end{proposition}

Our principal application of such diffeomorphism corks is the construction of exotic $\RR^4$'s equipped with exotic diffeomorphisms. The first examples of exotic $\RR^4$'s with nontrivial mapping class groups were due to Gompf \cite[Theorem 1.1]{gompf2018group}, although his examples are not compactly supported. The first example with $\MCG_c \neq 1$ was given in \cite[Theorem C]{konno2024four}. We extend this by constructing examples with arbitrarily large rank:

\begin{theorem}\label{thm:intro}
  For each $n \geq 1$, there exist uncountably many exotic copies of $\RR^4$ denoted $\R_{n,t}$ such that $\rank(\MCG_c(\R_{n,t})) \geq n$. 
\end{theorem}

To prove this, we apply a classical construction of Perron, which takes as input a contractible submanifold of a closed, simply-connected $4$-manifold with simply-connected exterior and outputs an open topological $4$-ball \cite[Lemma 7.2]{perron1986pseudo}. Applying this construction to a diffeomorphism cork $C$ as constructed by Proposition~\ref{prop:intro} produces an exotic copy of $\RR^4$ labeled $\R_C$. $C$ itself smoothly embeds in $\R_C$, so the diffeomorphisms supported on $C$ can be extended by the identity to compactly-supported diffeomorphisms of $\R_C$ \cite[Theorem C]{konno2024four}.

That $\R_C$ is exotic is intimately related to the fact that $C$ supports exotic diffeomorphisms, by way of the following result which is known to experts (cf. \cite[Remark 8.6]{krushkal2024corks}):
\begin{lemma}\label{lem:intro}
  Exotic diffeomorphisms detected by the families Seiberg-Witten invariant cannot be isotoped to have support in $B^4$.
\end{lemma}
We provide a proof of this folklore result --- to our knowledge, the first to appear in the literature --- and its natural generalization to the $S^1$-equivariant families Bauer-Furuta invariant in Section~\ref{sec:gaugetheory}. Given this technical result, if $\R_C$ were diffeomorphic to $\RR^4$, we should be able to isotope our exotic diffeomorphisms into a smooth $4$-ball, which provides the contradiction (this is spelled out in the proof of Theorem~\ref{thm:intro}). 

The exotic diffeomorphisms we use in the proof of Theorem~\ref{thm:intro} are due to Ruberman, who goes further and gives examples of closed, simply-connected $4$-manifolds supporting infinitely many independent exotic diffeomorphisms. However, we are only able to produce diffeomorphism corks (and subsequently exotic $\RR^4$'s) for finite subsets of this infinite collection. The following question is then natural:

\begin{question}
  Is there an exotic copy $\R$ of $\RR^4$ such that $\MCG_c(\R)$ is of infinite rank?
\end{question}

Gompf's examples, referenced above, include $\R$ such that $\MCG(\R)$ are of infinite rank. However, these diffeomorphisms are not compactly supported. 

In fact, our approach cannot be naïvely extended to the infinite rank case by a non-localization result of Konno-Mallick \cite{konno2024localizinggroupsexoticdiffeomorphisms}, which implies that none of Ruberman's examples of infinite families of exotic diffeomorphisms (or, indeed, all known infinitely generated groups of exotic diffeomorphisms of $4$-manifolds detected by families Seiberg-Witten theory) may be localized to even locally flatly embedded rational homology balls. 

\subsection{Notation and Conventions} We use $\homeo$ to denote homeomorphisms and $\diffeo$ to denote diffeomorphisms. $\MCG$ refers to the orientation preserving mapping class group, with decorations thereof (e.g. $\MCG_c$, $\MCG_\partial$). For concision, we adopt the standard abuse of notation that $A \setminus B$ really means $A \setminus \Int(B)$. 

\subsection*{Acknowledgements} The author thanks his advisor, Lisa Piccirillo, for helpful conversations and feedback. Thanks also to Imogen Montague for many useful discussions on gauge theory, and to Mark Powell and Anubhav Mukherjee for suggesting this project. This work was done with partial support by the Simons Collaboration grant `New Structures in Low-Dimensional Topology'.

\bigskip

\section{Diffeomorphism Corks}\label{sec:diffcorks}

\smallskip

\subsection{Background}
To modify the construction of diffeomorphism corks in service of Proposition~\ref{prop:intro}, we provide here a \emph{précis} of the relevant pseudoisotopy theory. Throughout, for reference, we refer to \cite[\S 2.1]{krushkal2024corks}. 

Morally, pseudoisotopies are to exotic diffeomorphisms as $h$-cobordisms are to exotic smooth structures (in dimension $4$), with Cerf theory replacing Morse theory. The construction of a sphere system from an $h$-cobordism (en route to the construction of a cork) has an analogue in the pseudoisotopy setting: the \emph{Quinn core}. The Quinn core and its properties will be the object of most of our discussion in this subsection.

Let $X$ be a closed simply-connected smooth $4$-manifold, and $F:X \times I \xrightarrow{\sim} X \times I$ a smooth pseudoisotopy i.e. a diffeomorphism of $X \times I$ restricting to the identity on $X \times \{0\}$. By Cerf theory, the data of $F$ is equivalent to a certain one-parameter family of generalized Morse functions $q_t:X \times I \rightarrow I$, called a \emph{Cerf family} for $F$. The endpoint functions $q_0(x,s) = s$ and $q_1(x,s) = q_0 \circ F^{-1}(x,s)$ have no critical points, so they correspond to the trivial handle decomposition of $X \times I$. 

Generalized Morse functions have isolated degenerate critical points, corresponding to births and deaths of cancelling pairs of handles. They are singular, but in the simplest possible way, as they lie in the codimension one stratum of Cerf's stratification of $C^\infty(X)$ \cite[Introduction, Section 4]{cerf1970stratification}. The input of Cerf theory here is that a generic path between $q_0$ and $q_1$ has at most these codimension one singularities, and if there are no such singularities along the path, then the pseudoisotopy $F$ is an isotopy.

We may impose more structure on $q_t$ (cf. \cite[\S 3.1]{krushkal2024corks}), in particular, we may arrange for all births to happen in the interval $t \in (\frac{1}{4}-\delta,\frac{1}{4})$ and all deaths in the interval $t \in (\frac{3}{4}, \frac{3}{4}+\delta)$. Thus, if there are $n$ births and deaths, for $t \in [\frac{1}{4}, \frac{3}{4}]$, $q_t$ gives $X \times I$ the structure of an $h$-cobordism with $n$ cancelling $2/3$ pairs. Note that $q_t^{-1}(\frac{1}{2}) := M_t \cong X \#_n S^2 \times S^2$, where the attaching ($A_i^t$) and belt ($B_i^t$) spheres of (respectively) the $3$ and $2$-handles intersect transversely. The number $n$ of births and deaths is referred to as the number of \emph{eyes} in a pseudoisotopy.

Just as the number of cancelling $2/3$-handle pairs in an $h$-cobordism determines the number of stabilizations needed to make an exotic pair stably diffeomorphic, the number of eyes in a pseudoisotopy determines the number of stabilizations needed to make an exotic diffeomorphism stably isotopic to the identity:

\begin{theorem}[{\cite[Theorems 2.5, 2.6]{krushkal2024corks}}]
  Let $f:X \rightarrow X$ be a diffeomorphism of a closed, simply-connected $4$-manifold. Then, there exists $n \geq 0$ such that the following are equivalent:
  \begin{enumerate}
    \item $f$ is topologically isotopic to $\id_X$
    \item $f$ is smoothly $n$-eyed pseudoisotopic to $\id_X$ (i.e., there exists a pseudoisotopy $F:X \times I \rightarrow X \times I$ such that $F|_{X \times \{1\}} = f$)
    \item $f$ is smoothly $n$-stably isotopic to $\id_X$
    \item $f$ is topologically pseudoisotopic to $\id_X$
  \end{enumerate}
\end{theorem}
We will restrict our attention to $1$-eyed pseudoisotopies where there is a single cancelling pair corresponding to a pair of spheres $A^t$, $B^t$, as this is the setting in which exotic diffeomorphisms can be localized to diffeomorphism corks.

Returning to the $q_t$, which describe a one-parameter family of trivial $h$-cobordisms, we want to understand how the description of the product structure on $X \times I$ given by $q_t$ changes over time. For $t \in [\frac{1}{4}, \frac{3}{4}]$, the $q_t$ taken together form a homotopy from $q_\frac{1}{4}$ to $q_\frac{3}{4}$. At $t=\frac{1}{4}$ and $t=\frac{3}{4}$, $A^t$ and $B^t$ intersect transversely in a single point. In the interval between them, the regular homotopy of $A^t \cup B^t$ (which restricts to an isotopy of $A^t$ and $B^t$ taken individually) introduces excess intersection between the spheres via finger moves, and removes them via Whitney moves. After a \emph{deformation} of the Cerf family (in the sense of \cite[\S~4.2]{quinn1986isotopy}, an isotopy of the pseudoisotopy), we may assume that the finger moves are all performed at $t= \frac{3}{8}$ along embedded finger discs $V_i$ and all the Whitney moves are performed at $t = \frac{5}{8}$ along embedded Whitney discs $W_j$. We will denote by $V,W$ the collection of finger and Whitney discs, and by $A,B$ the spheres $A^\frac{1}{2},B^\frac{1}{2}$ in $M := M_\frac{1}{2}$, the \emph{middle-middle level}. The finger and Whitney discs may intersect each other in their interiors, but there are no intersections between finger discs or between Whitney discs.

Then the \emph{Quinn core} of the Cerf family is defined as \[Q = \nu(A \cup B \cup V \cup W) \subseteq M\] We may then build a sub-$h$-cobordism $P \subseteq X \times I$ by taking $Q \times [-\epsilon,\epsilon]$ and attaching $3$-handles to either boundary component: to $Q \times \{\epsilon\}$, we attach them along the $A$ spheres, and to $Q \times \{-\epsilon\}$ (upside down) along the $B$ spheres. Setting $Q_i := P \cap (X \times \{i\}) \subseteq X$ for $i \in \{0,1\}$, by \cite[Lemma 3.2]{krushkal2024corks}, $F$ may be isotoped so that it is the identity on $(X \setminus Q_0) \times I$. Thus, the Quinn core isolates the nontriviality of the pseudoisotopy $F$. However, as we will see below, the Quinn core is not contractible.

We can put a handle structure on the Quinn core, coming from the $h$-cobordism structure on $X \times I$ induced by $q_\frac{1}{2}$. As above, we restrict our discussion to the one-eyed case.

\begin{enumerate}
  \item The spheres $A$ and $B$ are equipped with natural basepoints (labeled $*_A$ and $*_B$), namely, their original point of intersection with each other at $t=\frac{1}{4}$. By construction, $*_A = *_B$ so we take this point to be the basepoint of $Q$ (generically referred to as $*$), and take a neighborhood of $*$ in $M_\frac{1}{2}$ to be the $0$-handle for our decomposition.

  \item On $A$ and $B$, run disjoint arcs from $*$ to each point of \emph{excess} intersection with the other sphere. For each such point of intersection, there is a loop based at $*$ given by following the arc to it along $A$ and backwards along $B$. A neighborhood of this loop gives rise a $1$-handle attached to our $0$-handle. This adds $(|A \pitchfork B| - 1)$ $1$-handles to our handle decomposition.

  \item Each of $A,B$ minus a neighborhood of this tree of arcs (from the given sphere's basepoint to each excess intersection with the other sphere) is a disc, which we can thicken up to obtain $0$-framed $2$-handles labeled $H_A$ and $H_B$. The attaching circles for $H_A$ and $H_B$ are $0$-framed unknots which clasp each other at each point of intersection between $A$ and $B$. Since $A$ and $B$ intersect algebraically once, the attaching circles for $H_A$ and $H_B$ link (algebraically) once as well.

  \item Finally, the finger and Whitney discs can similarly be thickened up to obtain $0$-framed $2$-handles. Intersections between the finger and Whitney discs result, as above, in clasping of the corresponding attaching circles inside a dotted circle. This completes our handle decomposition for $Q$.
\end{enumerate}
  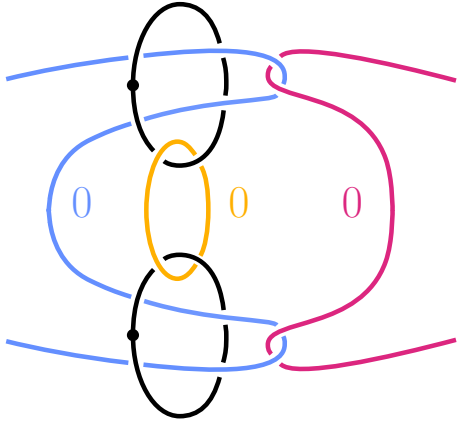
\begin{figure}[h]
    \centering
    \LARGE
    \def\svgwidth{.4\columnwidth}
\begingroup%
  \makeatletter%
  \providecommand\color[2][]{%
    \errmessage{(Inkscape) Color is used for the text in Inkscape, but the package 'color.sty' is not loaded}%
    \renewcommand\color[2][]{}%
  }%
  \providecommand\transparent[1]{%
    \errmessage{(Inkscape) Transparency is used (non-zero) for the text in Inkscape, but the package 'transparent.sty' is not loaded}%
    \renewcommand\transparent[1]{}%
  }%
  \providecommand\rotatebox[2]{#2}%
  \newcommand*\fsize{\dimexpr\f@size pt\relax}%
  \newcommand*\lineheight[1]{\fontsize{\fsize}{#1\fsize}\selectfont}%
  \ifx\svgwidth\undefined%
    \setlength{\unitlength}{196.28958875bp}%
    \ifx\svgscale\undefined%
      \relax%
    \else%
      \setlength{\unitlength}{\unitlength * \real{\svgscale}}%
    \fi%
  \else%
    \setlength{\unitlength}{\svgwidth}%
  \fi%
  \global\let\svgwidth\undefined%
  \global\let\svgscale\undefined%
  \makeatother%
  \begin{picture}(1,0.92452402)%
    \lineheight{1}%
    \setlength\tabcolsep{0pt}%
    \put(0,0){\includegraphics[width=\unitlength,page=1]{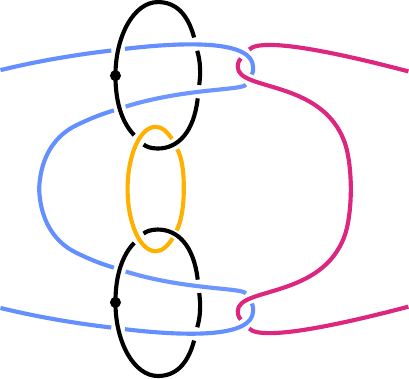}}%
    \put(0.49449283,0.44715003){\color[rgb]{1,0.69019608,0}\makebox(0,0)[lt]{\lineheight{1.25}\smash{\begin{tabular}[t]{l}$0$\end{tabular}}}}%
    \put(0.14599431,0.44715003){\color[rgb]{0.39215686,0.56078431,1}\makebox(0,0)[lt]{\lineheight{1.25}\smash{\begin{tabular}[t]{l}$0$\end{tabular}}}}%
    \put(0.7458736,0.44715003){\color[rgb]{0.8627451,0.14901961,0.49803922}\makebox(0,0)[lt]{\lineheight{1.25}\smash{\begin{tabular}[t]{l}$0$\end{tabular}}}}%
  \end{picture}%
\endgroup%

    \caption{The neighborhood of a pair of cancelling intersection points in a handle diagram (cf. \cite[Figure 3]{krushkal2024corks}). $\color{IBMblue}H_A$ and $\color{IBMred} H_B$ clasp each other through a dotted circle where the spheres $A$ and $B$ intersect. These two intersections are paired by a disc represented by the \textcolor{IBMyellow}{$0$-framed $2$-handle}, which may be either a finger or Whitney disc. Because there is no $1$-handle associated to the original point of intersection between $A$ and $B$, neither of the above clasps corresponds to the original point of intersection between $\color{IBMblue}A$ and $\color{IBMred}B$. The simplicity of this sub-diagram comes possibly at the expense of complicating the $2$-handle attaching curves for (respectively) the Whitney or finger discs.}
    \label{fig:handles}
\end{figure}
As Figure~\ref{fig:handles} may suggest, any complexity in a one-eyed Quinn core comes from intersections between the finger and Whitney discs. This intuition can be made precise as follows:

\begin{proposition}[{\cite[Corollary 3.8]{krushkal2024corks}}]\label{prop:homotopytype}
  The Quinn core of a one-eyed pseudoisotopy is homotopy equivalent to $S^2 \vee S^2 \vee \bigvee_1^m S^1$ where $m$ is the total number of interior intersections between the finger and Whitney discs.
\end{proposition}

In fact, the above statement can be deduced from the handle decomposition given above, together with the following final bit of structure we wish to impose on the Quinn core (see Remark~\ref{rem:htpytype}):
\begin{proposition}[{\cite[Lemma 3.7]{krushkal2024corks}}]\label{prop:arccondition}
  Let $Q$ be the Quinn core for a one-eyed pseudoisotopy. There exists a deformation of the Cerf data associated to this pseudoisotopy so that \emph{Quinn's arc condition} is satisfied, i.e., $(\partial V \cup \partial W) \cap A$ and $(\partial V \cup \partial W) \cap B$ are embedded arcs between the original point of intersection $A^\frac{1}{4} \pitchfork B^\frac{1}{4}$ and the final point of intersection $A^\frac{3}{4} \pitchfork B^\frac{3}{4}$.
\end{proposition}

\smallskip
\subsection{$\pi_1$-negligibility}

Having set the scene, we may outline our strategy for proving Proposition~\ref{prop:intro}. First, we consider the Quinn core of a single diffeomorphism, and deform its parent pseudoisotopy so that $\pi_1(M\setminus Q) = 1$ in the middle-middle level. We then apply a result of Melvin-Schwartz (quoted below as Lemma~\ref{lem:encasement}) to construct $R \supseteq Q$ in $M$ with $\pi_1(R) =1$ and $\pi_1(M\setminus R) = 1$. Analogous to the classical cork theorem, surgery on a $5$-dimensional handle attachment sphere will yield our diffeomorphism cork $C \subseteq X$.

For a finite collection $f_1, \cdots, f_n$ of diffeomorphisms, we build diffeomorphism corks $C_1,\cdots, C_n$ for each one as above, and consider their union in $X$. By another result of Melvin-Schwartz (quoted below as Lemma~\ref{lem:finger}), the $C_i$ may be repositioned by isotopies to maintain simple-connectivity of their collective exterior. Transverse intersections between the $C_i$ generate a free fundamental group for $\cup_{i=1}^n C_i$ which we may fill by another application of Lemma~\ref{lem:encasement}.

\begin{proposition}\label{prop:quinncore}
  After a suitable deformation of a $1$-eyed pseudoisotopy $F:X \times I \rightarrow X \times I$, the resulting Quinn core $Q$ in the middle-middle level $M$ can be taken to be $\pi_1$-negligible, i.e., $\pi_1(M \setminus Q) =1$.
\end{proposition}
\begin{proof}
  The proof will proceed in the following three stages:
  \begin{enumerate}
    \item First, we show that \[\pi_1(M \setminus \nu(A \cup B \cup V)) = \pi_1(M \setminus \nu(A \cup B \cup W)) = 1\]
    \item Moreover, $H_1(M\setminus Q) = 0$ (this is true without modification of $Q$).
    \item The above conditions imply that $\pi_1(M \setminus Q)$ is normally generated by meridians to the Whitney discs, and, separately, meridians to the finger discs. This allows us to kill $\pi_1(M\setminus Q)$ by finger moves between $V$ and $W$ (this is the only stage where we modify $Q$).
  \end{enumerate}

  We denote by $A',B'$ the attaching and belt spheres (respectively) of the unique pair of cancelling $3$ and $2$-handles shortly after their birth (i.e. $A' = A^\frac{1}{4}$). In particular, $A'$ and $B'$ intersect geometrically once. 

  First, note that $\pi_1(M \setminus \nu(A' \cup B')) = 1$ since $A'$ and $B'$ are geometrically dual; explicitly, the Seifert-van Kampen theorem tells us that $\pi_1(M \setminus \nu(A' \cup B'))$ is normally generated by meridians to $A'$ and $B'$, and since $A'$ and $B'$ are geometrically dual, $A'$ punctured at its intersection with $B'$ gives a nullhomotopy of $\mu_{B'}$ and vice versa. Pushoffs of $A'$ and $B'$ therefore give nullhomotopies of $\mu_{B'}$ and $\mu_{A'}$ (respectively) in the complement.

  Next, we can show that \[\pi_1(M \setminus \nu(A \cup B \cup V)) = \pi_1(M \setminus \nu(A \cup B \cup W)) = \pi_1(M \setminus \nu(A' \cup B')) = 1\]

  This follows from the fact that the complement of $\nu(A \cup B \cup V)$ is homotopy equivalent to the complement of $\nu(A' \cup B' \cup \upsilon)$ where $\upsilon$ is a collection of finger arcs guiding the finger moves. Explicitly, $\nu(A\cup B \cup V)$ is isotopic to $\nu(A' \cup B' \cup \upsilon)$ via a contraction of the disc which makes the surfaces tangent, and then pushing them off each other along a finger arc, so $\pi_1(M \setminus \nu(A \cup B \cup V)) = \pi_1(M \setminus \nu(A' \cup B' \cup \upsilon))$.

  Then, $\pi_1(M \setminus \nu(A' \cup B' \cup \upsilon)) = \pi_1(M \setminus \nu(A' \cup B'))$ since the finger arcs are $\pi_1$-negligible for dimensional reasons, or, equivalently, by noting that $M \setminus \nu(A' \cup B')$ is diffeomorphic to $M \setminus \nu(A' \cup B' \cup \upsilon) \cup 3\text{-handle}$. 

  The above argument with $V$ replaced by $W$ is identical with the caveat that contracting the Whitney discs rather than the finger discs results in the spheres $A^\frac{3}{4}$ and $B^\frac{3}{4}$, which intersect at a different point than $A'$ and $B'$, but are still geometrically dual, so the claim follows \emph{mutatis mutandis}. 

  Moreover, $H_1(M \setminus Q) = 0$. To see this, note that the $2$ and $3$-handle of the $h$-cobordism $X \times I$ (with the $h$-cobordism structure given by the Morse function $q_\frac{1}{2}$) are attached to $Q \times [-\epsilon,\epsilon]$, resulting in the sub-$h$-cobordism $Q_0 \times I \subseteq X \times I$. Therefore, the sub-$h$-cobordism above $X \setminus Q_0$ is \emph{canonically} a cylinder $(X \setminus Q_0) \times I \subseteq X \times I$ as it has no handles.

  Thus, $X \setminus Q_0 \diffeo M \setminus Q$ so $H_1(M \setminus Q) = H_1(X \setminus Q_0)$. We may apply the long exact sequence of the pair $(X,X\setminus Q_0)$ to obtain \begin{equation} \cdots \rightarrow H_2(X,X \setminus Q_0) \rightarrow H_1(X \setminus Q_0) \rightarrow H_1(X) = 0 \label{eq:LES}\end{equation} By excision, Poincaré-Lefschetz, and universal coefficients we know that \[H_2(X,X \setminus Q_0) \cong H_2(Q_0,\partial Q_0) \cong H^2(Q_0) \cong H_2(Q_0) \] In the last step, we know that $H_1(Q_0)$ has no torsion since, by Proposition~\ref{prop:homotopytype}, the homotopy type of $Q$ is a wedge of $1$ and $2$-spheres. $Q_0$ differs from $Q$ only where the $5$-dimensional handle attachment has the effect of dotting $H_B$. This $1$-handle (in $Q_0$) homologically cancels $H_A$ in $Q_0$, so $H_1(Q_0) = H_1(Q)$ has no torsion. 

  We also know that $H_2(Q_0) = 0$ since, by the aforementioned lemma, $H_2(Q) = \ZZ^2$, and by surgering the $B$-sphere, we kill $A$ in homology. Then (\ref{eq:LES}) implies that $H_1(X \setminus Q_0) = 0$. 

  Since $Q = \nu(A \cup B \cup V \cup W)$, the fact that $\pi_1(M \setminus \nu(A \cup B \cup V)) = 1$ implies that $\pi_1(M \setminus Q)$ is normally generated by meridians $\mu_W$ to the Whitney discs. Similarly, $\pi_1(M \setminus \nu(A \cup B \cup W)) = 1$ implies that $\pi_1(M \setminus Q)$ is normally generated by meridians $\mu_V$ to the finger discs. 

  Our goal is to trivialize each $\mu_V$ by modifying $Q$. Since $H_1(M \setminus Q) =0$, $\pi_1(M \setminus Q)/ \llangle \mu_V, \mu_W \rrangle = 1$ so we can write each $\mu_V$ as \[\mu_V = \prod_{i=1}^n [\mu_V^{v_i},\mu_W^{w_i}]^{\pm 1}\] where $\mu_V^{v_i} := v_i^{-1}\mu_V v_i$. The commutators on the right hand side above can be killed by finger moves between the appropriate finger and Whitney discs (see \cite{CASSONA1986Tlon} or \cite[\S 15.1]{behrens2021disc}), guided along arcs specified by the $v_i$ and $w_i$. Such finger moves constitute an allowable \emph{deformation} of the data of the pseudoisotopy per \cite[\S 4.2]{quinn1986isotopy}, which give rise to a modified Quinn core $Q$ satisfying $\pi_1(M \setminus Q) = 1$.
\end{proof}
\begin{remark}
  Note that we need that $\pi_1(M \setminus Q)$ is normally generated by both $\mu_V$ and $\mu_W$ to write $\mu_V$ as a product of commutators as above, since we want to deform $Q$ by finger moves between the finger moves and Whitney discs. If we used only that $\pi_1(M \setminus Q) = \llangle \mu_V \rrangle$, then killing $\mu_V$ by finger moves as above would possibly entail self-finger moves of a disc, which is not a valid deformation of the data.
\end{remark}

To build a contractible handlebody $C$ containing $Q$ while maintaining simple-connectivity of its exterior, we will apply a result of Melvin-Schwartz. Morally, their result is a generic version of the arguments in the proof of the classical cork theorem \cite{Curtis1996}, where extra handles are added to produce a cork from the sphere system associated to an $h$-cobordism, and further handles are added and slid to make the exterior of this cork simply-connected. To precisely state their result, we need the following definition:

\begin{definition}[AC Manifolds {\cite[Definition 1.1]{melvin2021higher}}]
  Let $X$ be a $4$-manifold with only $0$, $1$, and $2$-handles. We say that a handle decomposition is $\emph{AC}$ (for Andrews-Curtis) if some subset of the $2$-handles homotopically cancel some subset of the $1$-handles, and the remaining $2$-handles are attached homotopically trivially. A $4$-manifold is AC if it has an AC handle decomposition.

  $X$ is AC of Type I if all $2$-handles homotopically cancel some $1$-handle. $X$ is AC of Type II if all $1$-handles are homotopically cancelled by some $2$-handle. 
\end{definition}
\begin{lemma}[{Encasement Lemma \cite[Lemma 2.2]{melvin2021higher}}]\label{lem:encasement}
  Let $X$ be a closed simply-connected $4$-manifold, and $U$ any AC manifold lying in $X$. Then there exists an AC manifold $V$ of type II with $U \subseteq V \subseteq X$ whose AC structure extends that of $U$, inducing a canonical isomorphism $H_2(U) \cong H_2(V)$ sending generators to generators. Moreover, if $X \setminus U$ is simply-connected, then $V$ can be chosen so that $X \setminus V$ is AC as well.
\end{lemma}

\begin{remark}
  Any AC manifold is homotopy equivalent to a wedge of circles and spheres. Type I AC manifolds are homotopy equivalent to a bouquet of circles, and type II AC manifolds to a bouquet of spheres. 

  The notion of an AC manifold is implicit in the proof of the original cork theorem; in particular, a major step in the proof is to turn a handle decomposition built from the $h$-cobordism into an AC handle decomposition. The key idea is that in an AC manifold, $2$-handles are neatly sorted into two roles: cancelling a $1$-handle, or generating $H_2$, whereas, in nature, $2$-handles are not so cleanly delineated. 

  For example, the standard handle decomposition for the Akbulut cork is an AC structure, giving rise to the fundamental group presentation $\pi_1(A) = \langle x | x x^{-1} x \rangle = 1$. Importantly, the relator defined by the $2$-handle is homotopic to $x$, but is not literally the word $x$, which would correspond to a cancelling $2$-handle.
\end{remark}
\begin{proposition}\label{prop:QisAC}
  The Quinn core of a one-eyed pseudoisotopy has an AC structure.
\end{proposition}
\begin{proof}
  The handle decomposition given in the discussion preceding Figure~\ref{fig:handles} is not quite AC, but can easily be made AC by a sequence of handle slides. First, we assume that Quinn's arc condition (Lemma~\ref{prop:arccondition}) holds, so we may label our intersection points $v_1,\cdots,v_{2k+1}$ and our finger and Whitney discs as $V_1,\cdots,V_k$ and $W_1,\cdots, W_k$. This labeling has the useful property that $W_i$ pairs the intersections $v_{2i}$ and $v_{2i-1}$, and $V_i$ pairs the intersections $v_{2i+1}$ and $v_{2i}$. 

  Let $g_i$ be the $1$-handle generator corresponding to the intersection $v_i$ for $i \geq 2$, and $h_1,\cdots,h_m$ the $1$-handle generators corresponding to intersections between the finger and Whitney discs. Note that there is no $1$-handle corresponding to $v_1$ in our handle decomposition by choice of basepoint for $A$ and $B$. Then, the presentation of $\pi_1(Q)$ induced by our handle decomposition is calculated in the proof of \cite[Lemma 3.5]{krushkal2024corks} as \[\pi_1(Q) = \langle g_2,\cdots,g_{2k+1},h_1,\cdots,h_m| g_2, g_2^{-1}g_3,\cdots,g_{2k}^{-1}g_{2k+1}, 1, 1 \rangle\] where the two trivial relators correspond to $H_A$ and $H_B$. Note that we can slide $2$-handles to concatenate their relators, so the second relator becomes $g_3$ after such a move, and, inductively, all remaining relators become of the form $g_i$ for some $i$. The resulting presentation of $\pi_1(Q)$ corresponds to an AC handle decomposition for $Q$.
\end{proof}
\begin{remark}\label{rem:htpytype}
  The proof that $Q$ can be made AC gives an alternate proof of Proposition~\ref{prop:homotopytype}, since the homotopically trivial relators corresponding to $H_A$ and $H_B$ give $S^2$ wedge factors, and the remaining free $\pi_1$ generators $h_i$ correspond to the $S^1$ wedge factors.
\end{remark}

We now recall the existence theorem for diffeomorphism corks, as stated and proven in \cite{krushkal2024corks}:
\begin{theorem}[{\cite[Theorem 4.1]{krushkal2024corks}}]\label{thm:4.1}
  For $F:X \times I \rightarrow X \times I$ a $1$-eyed pseudoisotopy, there exists a smooth, compact, contractible, codimension zero submanifold $C \subseteq X$ and a smooth isotopy fixing $X \times \{0\}$ from $F$ to a pseudoisotopy supported on $C \times I$.
\end{theorem}
This can be improved as follows:
\begin{proposition}\label{prop:cork}
  $C$ may be taken to be AC, and satisfying $\pi_1(X \setminus C) = 1$. 
\end{proposition}
\begin{proof}  
  Since $Q$ is AC by Proposition~\ref{prop:QisAC}, we may apply Lemma~\ref{lem:encasement} to the Quinn core $Q \subseteq M$ produced by Proposition~\ref{prop:quinncore} to obtain $R$ satisfying
  \begin{enumerate}
     \item $Q \subseteq R$
     \item $R$ is AC of type II, so $\pi_1(R) = 1$
     \item $H_2(R) \cong H_2(Q) \cong \ZZ^2$
     \item $M \setminus R$ is AC
   \end{enumerate}
  Let $U$ be given by $R \times [-\epsilon,\epsilon]$ together with the $5$-dimensional $3$-handles of the $h$-cobordism attached to $A \subseteq R \times \{\epsilon\}$ and (upside down) to $B \subseteq R \times \{-\epsilon\}$. These handles geometrically cancel the generators of $H_2(R) \cong H_2(Q) \cong \ZZ^2$ so $U$ is a contractible sub-$h$-cobordism of $X \times I$. 

  We can then define $C := U \cap (X \times \{0\})$; observe that $Q_0 \subseteq C$ by construction. Moreover, since $(X \times I) \setminus U$ contains no handles, it is a cylinder, so $X \setminus C \diffeo M \setminus R$, which implies that $H_1(M \setminus R) = H_1(X \setminus C)$. As we will show below, $C$ is contractible, so $H_1(X \setminus C) = 0$ by Mayer-Vietoris, hence $H_1(M \setminus R) = 0$. Since $M \setminus R$ is AC, and has vanishing first homology, it must be AC of type II, whence $\pi_1(M \setminus R) = \pi_1(X \setminus C) =  1$.

  $C$ is contractible and AC; stated differently, $C$ is AC of type II with $H_2(C) = 0$. For the latter claim, note that the handle diagram of $C$ derived from that of $R$ (see Figure~\ref{fig:handles} and preceding discussion). The only difference between the handle diagrams is that the $2$-handle attaching curve for $H_B$ is replaced \emph{in situ} with a dotted circle for a $1$-handle. This new $1$-handle cancels $H_A$ in homology, so $H_2(C) = 0$. 


  To see that $C$ is AC of type II, we need to study the effect of this new $1$-handle on the AC presentation of $R$. By Proposition~\ref{prop:QisAC}, we naturally obtain a presentation of $\pi_1(C)$ from that of $\pi_1(Q)$ with one new generator $t$ and where the relators $g_i^{-1}g_{i+1}$ becomes $g_i^{-1}t^{a_i}g_{i+1}t^{b_i}$. The relator $r_A$ corresponding to $H_A$ was trivial in $Q$, and in $C$ becomes a product in the $g_i$ and $t$ with total exponent $A \cdot B = 1$ in $t$; we may replace all instances of $g_{2k+1}$ in this word with $g_{2k}$ (along with some powers of $t$) using the above relations, and inductively it follows that $r_A =t$ in $\pi_1(C)$, and the new $1$-handle is homotopically cancelled by $H_A$. For any other $2$-handles that may have linked $H_B$ in $R$, the relator they spell in $C$ is modified by the presence of the new generator, but the AC condition is preserved since this new generator is itself nullhomotopic. $C$ is then AC of type II since $R$ is.

  We already know that $F$ may be isotoped to be supported on $Q_0 \subseteq C$ by \cite[Lemma 3.2]{krushkal2024corks}, so in particular, $F$ is supported in $C$.
\end{proof}

The following lemma of Melvin-Schwartz then suffices to extend the above result to finite collections of diffeomorphisms:
\begin{lemma}[{Finger Lemma \cite[Lemma 2.4]{melvin2021higher}}]\label{lem:finger}
  Any finite collection of contractible AC manifolds $C_1,\cdots,C_n$ embedded in a $4$-manifold $X$ can be repositioned by isotopies so that their union $\cup_{i=1}^n C_i$ is an AC manifold of type $I$ satisfying
  \begin{enumerate}
    \item The AC structure on each $C_i$ extends to an AC structure on $\cup_{i=1}^n C_i$
    \item $X \setminus \nu(\cup_{i=1}^n C_i)$ is simply-connected provided that each $X \setminus C_i$ is simply connected
  \end{enumerate}
\end{lemma}
\begin{remark}\label{rem:dblptloops}
  In general, $\cup_{i=1}^n C_i$ has a nontrivial fundamental group arising from intersections between the $C_i$. We may assume that such intersections take place in the core discs of the $2$-handles by transversality. Such intersections (assuming, as we will below, that $\cup_{i=1}^n C_i$ is connected) give rise to a free fundamental group generated by \emph{double point loops} \cite[\S 1.2]{freedman2014topology} through the intersections; such a loop is given by following a path to the point of intersection along one sheet, and returning to the basepoint along the other sheet. 

  The generators for this free fundamental group can alternatively be seen implicitly in Figure~\ref{fig:handles}, where the intersection of two essential surfaces can be locally modeled in a handle diagram by clasping the corresponding handles inside a dotted circle. This $1$-handle (with no attached finger/Whitney disc) is the free generator corresponding to the intersection. 
\end{remark}

As above and for clarity, we recall the analogous theorem from \cite{krushkal2024corks}: 
\begin{theorem}[{\cite[Theorem 4.2]{krushkal2024corks}}]\label{thm:4.2}
  For any finite collection $\{f_i\}_{i=1}^n$ of self-diffeomorphisms of a closed, simply-connected, smooth $4$-manifold $X$, each $1$-stably isotopic to the identity, there exists a smooth, compact, contractible, codimension zero submanifold $C \subset X$ and a collection $\{f_i'\}_{i=1}^n$ such that $f_i'$ is smoothly isotopic to $f_i$, and the $f_i'$ are all supported on $C$. 
\end{theorem}
This can be improved as follows:
\begin{lemma}\label{lem:finitecollections}
  In the above, we may take $C$ to be AC and satisfying $\pi_1(X \setminus C) = 1$.
\end{lemma}
\begin{proof}
  For each $f_i$ we may apply Proposition~\ref{prop:cork} to produce a contractible AC manifold $C_i$ with simply-connected exterior such that there exists $\tilde f_i$ isotopic to $f_i$ supported on $C_i$. Then, by Lemma~\ref{lem:finger}, the $C_i$ may be repositioned by isotopies so that $\cup_{i=1}^n C_i$ (which we take to be connected by finger moves between the core discs of the $2$-handles of the $C_i$ if necessary) is an AC manifold of type I with simply-connected exterior. We may isotope the $\tilde f_i$ along the isotopies repositioning the $C_i$ to produce the $f_i'$.

  Then, we may again apply Lemma~\ref{lem:encasement} to produce a contractible AC manifold $C \supseteq \cup_{i=1}^n C_i$. Since $H_1(X \setminus C) = 0$, and $X \setminus C$ is AC, $\pi_1(X \setminus C) = 1$ as desired.
\end{proof}

\bigskip

\section{Exotic $\RR^4$'s}\label{sec:exotica}

\smallskip

In this section, we prove Theorem~\ref{thm:intro} by constructing exotic copies $\R_n$ of $\RR^4$ containing diffeomorphism corks. The induced diffeomorphisms on the $\R_n$ are nontrivial (i.e., not isotopic to the identity), which establishes that the compactly-supported mapping class group $\MCG_c(\R_n)$ is nontrivial. 
\begin{definition}[$\MCG_c$]
  The \emph{compactly supported mapping class group} $\MCG_c(X)$ of a manifold $X$ is the group of (orientation preserving) compactly supported diffeomorphisms of $X$ modulo smooth compactly-supported isotopy.
\end{definition}
For example, $\MCG_c(\RR^4) \cong \MCG_\partial(B^4) \cong \MCG(S^4)$. 

For our construction of the $\R_n$, we need the following lemma due originally to Andrews-Curtis \cite[Corollary 1]{andrews1965free}:

\begin{lemma}\label{lem:doublestos4}
  The double of a contractible AC manifold is $S^4$.
\end{lemma}
\begin{proof}
  The point is that $C \cup_\partial -C = \partial (C \times I)$ and $C \times I \diffeo B^5$. This follows from the observation that homotopy implies isotopy for curves in $4$-manifolds, so the homotopic cancellation between the $1$ and $2$-handles of $C$ (which comes from the fact that $C$ is AC) may be upgraded to geometric cancellation of the corresponding handles in $C \times I$, whence $C \times I \cong D^5$.
\end{proof}
\begin{theorem}\label{thm:main}
  For each $n \geq 1$, there exists an exotic copy of $\RR^4$ denoted $\R_n$ and a collection $\{g_1,\cdots,g_n\}$ of compactly supported diffeomorphisms of $\R_n$ which are independent in the sense that $\rank(\MCG_c(\R_n)) \geq n$. 
\end{theorem}
\begin{proof}
  Our examples are derived from Ruberman's \cite{ruberman1999polynomial}. For any $m \geq 2$, he constructs an infinite family of exotic diffeomorphisms of \[Z_m := \begin{cases} \#_{2m} \CC \PP^2 \#_{10m+1} \overline{\CC\PP^2} & m \text{ odd} \\ \#_{2m}\CC\PP^2 \#_{10m+2} \overline{\CC\PP^2} & m \text{ even} \end{cases}\] In particular, he constructed a group homomorphism $\DD: \MCG^T(Z_n) \rightarrow \RR[\![H_2(Z_n)^\ast]\!]$ from the Torelli subgroup (of diffeomorphisms acting trivially on homology) to a power series ring, and shows that the diffeomorphisms are linearly independent in the image of this map. Note that, by \cite[Theorem C]{auckly2015stable}, \cite[Example 7.4]{krushkal2024corks}, these diffeomorphisms are $1$-stably isotopic to the identity.

  For convenience, fix $m = 2$ with $Z := Z_2$, and consider $f_1,\cdots,f_n$ the first $n$ of Ruberman's family of diffeomorphisms of $Z$. By Lemma~\ref{lem:finitecollections}, the collection $\{f_1,\cdots,f_n\}$ can be isotoped to a collection $\{g_1,\cdots,g_n\}$ such that the $g_i$ are all supported on a smooth, compact, contractible, codimension zero submanifold $C_n$ of $Z$. By \cite[Theorem~1.6]{krushkal2024corks}, the $g_i$ generate a subgroup $D$ of $\MCG_\partial(C_n)$ which abelianizes to $\ZZ^n$, hence $\rank(D) = n$. 

  By \cite{boyer1986simply,stong1993simply}, $Z \setminus C_n$ is homeomorphic rel. boundary to $Z \# \overline{C_n}$ since both manifolds are simply-connected with the same intersection form and integer homology $3$-sphere boundary. Let $\rho':Z \setminus C_n \rightarrow Z \# \overline{C_n}$ be such a homeomorphism. Then, extending by the identity, we obtain a homeomorphism $\rho:(Z \setminus C_n) \cup_\partial C_n \rightarrow (Z \# \overline{C_n}) \cup_\partial C_n$. The domain of this map is $Z$, and the codomain is \[(Z \# \overline{C_n}) \cup_\partial C_n = Z \# (\overline{C_n} \cup_\partial C_n) = Z \# S^4\] where $\overline{C_n} \cup_\partial C_n \diffeo S^4$ by Proposition~\ref{lem:doublestos4}. 

  Set $\R_n \subseteq Z$ to be the open subset \[\R_n := \rho^{-1}(\Int(\overline{C_n} \cup_\partial C_n)^\circ)\] equipped with the smooth structure from the domain --- explicitly, the smooth atlas on $\R_n$ is given by restriction of the smooth atlas on $Z$, which is well-defined since $\R_n$ is open in $Z$. Now $\rho$ restricts to a homeomorphism from $\R_n$ to the standard $\RR^4$ since $\Int(\overline{C_n} \cup_\partial C_n)^\circ \diffeo \RR^4$, so $\R_n \homeo \RR^4$. 

  The key property of $\R_n$ is that the $g_i$ naturally act on it. Note that $\rho'$ is the identity on $\partial C_n$, and $\rho$ is the identity by construction on $C_n$, so $\rho^{-1}(C_n) \diffeo C_n$ is naturally embedded in $\R_n$. Thus, we may extend the $g_i$ by the identity outside of $C_n$ to obtain compactly-supported diffeomorphisms of $\R_n$. 

  Moreover, since $\R_n \hookrightarrow Z$, if $h \in \MCG_c(\R_n)$ is compactly-supported isotopic to the identity, the induced diffeomorphism of $Z$ is isotopic to the identity as well. It follows that the map $\MCG_c(\R_n) \rightarrow \MCG^T(Z)$ induced by inclusion is injective on the $g_i$, so we may define a composite homomorphism \[\MCG_c(\R_n) \xrightarrow{\iota_\ast} \MCG^T(Z) \xrightarrow{\DD} \RR[\![H_2(Z)^\ast]\!]\] such that the images of the $g_i$ are linearly independent, from which it follows that $\rank(\MCG_c(R_n)) \geq n$.

  It remains to show that $\R_n \centernot \diffeo \RR^4$. This will follow from Corollary~\ref{cor:FSWfolklore} below, which states that exotic diffeomorphisms detected by families Seiberg-Witten theory (which our $f_i$ are, per \cite[Corollary 1.3]{baraglia2020gluing}) cannot be isotoped to have support in a $4$-ball. If $\R_n \diffeo \RR^4$, then we can choose a sufficiently large $4$-ball containing $C_n \hookrightarrow \R_n$ and obtain a contradiction. 
\end{proof}
\begin{remark}
  Note that the above proof suffices to show that any finite collection of exotic diffeomorphisms that are $1$-stably isotopic to the identity and detected by the families Seiberg-Witten invariant induce a finite collection of compactly supported exotic diffeomorphisms (up to compactly supported isotopy) of some exotic $\RR^4$. However, the bound on the rank of $\MCG_c(-)$ made essential use of the homomorphism $\DD$ and does not automatically generalize.

  Note also that the above proof suffices to show that $\MCG_c(\R_n)/\MCG_\partial(B^4)$ has rank at least $n$, where $\MCG_\partial(B^4)$ is a well-defined (possibly trivial) subgroup of $\MCG_c(\R_n)$ by the Cerf-Palais disc theorem. Thus, $\R_n$ and $\RR^4$ are distinguished by $\MCG_c(-)/\MCG_\partial(B^4)$. 

  Given the inhomogeneity of all known exotic $\RR^4$'s, one might expect that $\MCG_c$ is often large, reflecting the large amount of varying ``local'' structure to play with. Perhaps counter to this intuition, it has been announced \cite[Theorem 1]{RayGompfOrson2026} that the Freedman-Taylor \cite{freedman1986universal} universal exotic $\RR^4$ has trivial compactly supported mapping class group, the first such known example (since $\MCG_c(\RR^4_\std) = \MCG(S^4)$).
\end{remark}

Following a strategy of Gompf, we may pass from the $\R_n$ constructed above to uncountable families thereof:
\begin{theorem}\label{thm:uncountable}
  The $\R_n$ in Theorem~\ref{thm:main} may be augmented to an uncountable family $\R_{n,t}$ of pairwise non-diffeomorphic exotic $\RR^4$'s with $\rank(\MCG_c(\R_{n,t})) \geq n$.
\end{theorem}
\begin{proof}
  The main ingredients in this improvement are an exotic $\RR^4$, $\R_L \subseteq \overline{\CC\PP^2}$ from \cite[\S 3]{gompf2018group}, and a family Donaldson invariant blowup formula \cite[Theorem 3.2]{ruberman1999polynomial}. The basic idea is that we end sum $\R_n$ and $\R_L$ to obtain an exotic $\RR^4$ contained in $Z_2 \# \overline{\CC\PP^2} = \#_4 \CC \PP^2 \#_{23}\overline{\CC\PP^2}$. We can then use the standard radial family technique to produce an uncountable family of exotic $\RR^4$'s which are pairwise nondiffeomorphic by Gompf's result, and in which our exotic diffeomorphisms $g_i$ remain exotic by the blowup formula.

  The key properties of $\R_L \subseteq \overline{\CC\PP^2}$ that we will require are that $\R_L$ is diffeomorphic at infinity to $V_L$, a simply-connected, spin, positive-definite $4$-manifold with nontrivial (and therefore non-diagonalizable) intersection form, and that $\R_L$ contains a compact subset $K_L$ that does not embed in any closed, positive-definite $4$-manifold (exoticity of $\R_L$ then follows by Donaldson's theorem). By diffeomorphic at infinity, we mean that there exist compact subsets $A \subseteq \R_L$ and $B \subseteq V_L$ and a diffeomorphism $\psi:\R_L \setminus A \xrightarrow{\sim} V_L \setminus B$; two such diffeomorphisms are equivalent if they agree outside some larger compact subset. 

  Let $\R' := \R_n \natural \R_L \subseteq Z_2 \# \overline{\CC\PP^2}$ where $\natural$ denotes the end sum (and where we choose the end-summing rays to avoid $C_n \subseteq \R_n$ and $K_L \subseteq \R_L$), and pick a homeomorphism $f:\R' \rightarrow \RR^4_\std$. Since $C_n$ and $K_L$ are both compact, there exists $r_0 > 0$ such that $f^{-1}(B_{r_0}(0))$ contains them, where $B_r(x)$ denotes the open ball centered at $x$ with radius $r$; set $\R_{n,t} := f^{-1}(B_t(0))$ for all $t \geq r_0$. 

  We now follow a nearly identical argument to that in \cite[\S 3]{gompf2018group}, which we outline here for convenience. The $\R_{n,t}$ are pairwise nondiffeomorphic; this follows from Gompf's observation that, given $K_L \subseteq \R_a \subseteq \R_b$ with the closure of $\R_a$ compact in $\R_b$, $\R_a$ and $\R_b$ have nondiffeomorphic end neighborhoods; otherwise, enlarging $K_L$ as necessary, we may pick a diffeomorphism of ends $\phi:\R_L \setminus A \xrightarrow{\sim} V_L \setminus B$ where $A \subseteq \Int(K_L)$ and where $H_2(B) = H_2(V_L)$ and form an open manifold with periodic ends by gluing infinitely many copies of the region between $\R_a$ and $\R_b$ (whose ends are diffeomorphic) to $(\R_L \setminus A) \cup_{\phi|_{\partial A}} B$. Such a manifold would contradict Taubes' theorem \cite[Theorem 1.4]{taubes1987gauge}, a diagonalization theorem for manifolds with periodic ends.

  Moreover, the $g_i$ survive as exotic diffeomorphisms to $\R_{n,t}$, since the form of the blowup formula easily shows that any relation among the $g_i \# \id_{\overline{\CC\PP^2}}$ induces one among the $g_i$. 
\end{proof}

\section{Families Bauer-Furuta Invariants}\label{sec:gaugetheory}
In this section, we prove Lemma~\ref{lem:intro}, i.e., we prove that diffeomorphisms detected by families Seiberg-Witten theory cannot be isotoped to be supported in $B^4$; this will follow from a connect-sum formula for the \emph{families Bauer-Furuta invariant}.

Following \cite[\S 2.2]{lin2023isotopy} we now outline the construction of the $S^1$-equivariant families Bauer–Furuta invariant of diffeomorphisms $\FBF^{S^1}(X,f)$ when $b_1(X) = 0$, and $f:X \rightarrow X$ is an orientation-preserving diffeomorphism acting trivially on $H_2(X)$ (and, in particular, preserving the spin$^c$ structure $\s$ on $X$). We will indicate along the way how the general spin$^c$ case differs from the spin case, which is more thoroughly treated in the literature. 

Consider the mapping torus $M_f$ of $f$ as an $X$-bundle over $S^1$ equipped with a family of fiberwise Riemannian metrics that are standard on the collars $N_i$. The automorphism group of a spin$^c$ structure on $X$ is the gauge group $\G$ of maps from $X$ to $S^1$, which splits as $\G = \G_0 \times S^1$ where $\G_0$ consists of based gauge transformations. Since $b_1(X) = 0$, and $S^1$ is a $K(\ZZ,1)$, $\G_0$ and $\G$ are connected (in general $\pi_0(\G_0) \cong H^1(X;\ZZ) = \ZZ^{b_1}$ with each component contractible). Thus, by the clutching construction, each spin$^c$ structure on $X$ lifts uniquely to a fiberwise spin$^c$ structure on $M_f$. In contrast, in the spin case (also when $b_1 = 0$) there are two fiberwise spin structures corresponding to the unique spin structure.

Note that Kronheimer-Mrowka \cite[Definition 3.2]{Kronheimer2020-jo} work over an arbitrary base $B$ rather than $S^1$ which necessitates the additional assumption that $\pi_1(B)$ acts trivially on $H_2(X)$; this assumption is subsumed under our requirement that $f$ is in the Torelli subgroup (i.e., acts trivially on $H_2(X)$). 

With these data, we can construct the fiberwise monopole map $\SW: \mathcal{W}^+ \rightarrow \mathcal{W}^-$ as a map of Hilbert bundles over $S^1$. $\mathcal{W}^\pm$ decomposes as $\mathcal{W}^\pm = \mathcal{V}^\pm \oplus \mathcal{U}^\pm$ where the $\mathcal{V}^\pm$ are Sobolev completions of $\Gamma(S^\pm)$, $\mathcal{U}^+$ is a Sobolev completion of $\Omega^1(X)$, and $\mathcal{U}^-$ is a Sobolev completion of $\Omega^2_+(X) \oplus \Omega^0(X)/\RR$. 

  The procedure for trivializing these bundles and applying finite-dimensional approximation is detailed in \cite[2.2.2]{lin2023isotopy}, and we may follow the same steps in the spin$^c$ setting, replacing the universe generated by the $\Pin(2)$-representations $\HH$, $\tilde \RR$, and $\RR$ by one generated by the $S^1$-representations $\CC$ and $\RR$. The finite-dimensional approximations to $\mathcal{V}^\pm$ are denoted by $V^\pm$ and similarly for $\mathcal{U}^\pm$ and $\mathcal{W}^\pm$, and we may canonically trivialize $U^+$ and $V^-$ using Kuiper's theorem and its equivariant generalization \cite{kuiper1965homotopy,segal1969equivariant}.

  The principal difference in the spin$^c$ case is in the trivialization of $V^+$; in particular, $S^\pm$ are $\U(2)$ bundles as opposed to $\Sp(1)$ bundles, so the finite-dimensional approximation to $\Gamma(S^\pm)$ is a $\U(N)$ bundle. In the spin setting, $\pi_0(\Sp(N)) = \pi_1(\Sp(N)) = 1$ implies that the bundle $V^+$ is trivial over $S^1$, and canonically so up to homotopy. In the spin$^c$ case, $\pi_1(\U(N)) = \ZZ$, which introduces a choice in the trivialization of $V^+$. Fortunately, \cite[Proposition 2]{szymik2010characteristic} implies that this choice does not affect the invariant we will ultimately assign to $f$. 

  Thus, taking the Thom space of the domain and projecting to a (compactified) fiber in the codomain yields a well-defined $S^1$-equivariant stable cohomotopy element \[[\sw] \in \{S^{\frac{c_1^2(\s) - \sigma(X)}{8}\CC} \wedge S^1_+, S^{b_2^+(X)\RR}\}^{S^1}\] where $S^1_+$ denotes $S^1$ together with a disjoint basepoint. To instead obtain an equivariant stable homotopy element, we define \[\FBF^{S^1}(X,f,\s) := p^\ast[\sw] \in \{S^{\RR + \frac{c_1^2(\s) - \sigma(X)}{8}\CC} , S^{b_2^+(X)\RR}\}^{S^1}\] where $p:S^{2\RR} \rightarrow S^\RR \wedge S^1_+$ is the Pontryagin-Thom collapsing map for $S^1 \hookrightarrow S^{2\RR}$. In practice, we will generally suppress the spin$^c$ structure $\s$ in our notation. 

We now turn to the setup for the connect-sum formula. Following the setup in \cite[\S 2]{bauer2004stableII}, let $X = \sqcup_{i=1}^n X_i$ where each $X_i$ is a closed, oriented, Riemannian $4$-manifold with $b_1 = 0$, with a spin$^c$-structure $\s_i$ on each $X_i$. Moreover, each $X_i$ splits as $X_i = X_i^- \cup_{S^3} X_i^+$, so for an even permutation $\tau \in A_n$, we define $X^\tau$ to be the manifold with connected components $X_i^\tau = X_i^- \cup_{S^3} X_{\tau(i)}^+$. The permutation $\tau$ is required to be even in Bauer's setup as the corresponding matrix needs to have determinant one for orientations to glue appropriately.

The orientations and spin$^c$ structures on the $X_i$ glue to orientations and spin$^c$ structures on the $X_i^\tau$. Suppose further that $f:X \rightarrow X$ is a diffeomorphism, compatible with the $\pm$-decomposition of each component $X_i$ in the sense that $f$ restricts to the identity on a bicollar $N_i := S^3 \times I$ of each separating $S^3$ giving rise to the decomposition $X_i = X_i^- \cup_{S^3} X_i^+$. Thus, there is a natural diffeomorphism $f^\tau:X^\tau \rightarrow X^\tau$ defined by gluing the restrictions of $f$ to $X_i^\pm$ along the identity on the $N_i$.

\begin{lemma}\label{lem:connectsum}
  \[\FBF^{S^1}(X,f) = \FBF^{S^1}(X^\tau,f^\tau)\] 
\end{lemma}

\begin{remark}
  The above lemma appears verbatim in \cite[Theorem 34]{xu2004bauer}; however, they assume that the orbit of $\s$ under $f^\ast$ is infinite, and their definition of $\FBF^{S^1}$ is not \emph{a priori} equivalent to our definition, which is now standard in the literature. Roughly, they consider the mapping cylinder of $f$ as opposed to the mapping torus, and their approach is more analogous to the definition of $\FSW$ in \cite{Ruberman_1998}. 
\end{remark}

\begin{proof}
  With minor modifications, the proof essentially follows the outline of \cite[Theorem 2.1]{bauer2004stableII}. We have supplied the data necessary to define the monopole map as a map of Hilbert bundles over $S^1$ for $(X,f)$ and $(X^\tau,f^\tau)$. Bauer's proof proceeds by constructing comparison maps $V$ between $\mathcal{W}^\pm$ and $(\mathcal{W}^\pm)^\tau$ then exhibiting a three-stage homotopy between $\SW$ (the monopole map for a manifold, not a family) and $V^{-1} \SW^\tau V$ satisfying certain boundedness criteria, from which one may deduce that the finite-dimensional approximations of $\SW$ and $\SW^\tau$ agree.

  As noted in \cite{Kronheimer2020-jo,lin2023isotopy}, this proof readily carries over to the families setting, since the homotopies can be applied fiberwise over $S^1$, with the estimates required for the boundedness criteria extended uniformly over the compact base of the mapping torus. These homotopies are $S^1$-equivariant by construction, so the finite-dimensional approximations of $\SW$ and $\SW^\tau$ agree and the result follows.
\end{proof}
\begin{corollary}\label{cor:connectsumidentity}
    \[\FBF^{S^1}(X_1 \# X_2, f \# \id_{X_2}) = \FBF^{S^1}(X_1,f) \wedge \BF^{S^1}(X_2)\] with appropriately restricted spin$^c$ structures.
\end{corollary}
\begin{proof}
  Set \[X = (X_1^\circ \cup_{S^3} B^4) \sqcup (B^4 \cup_{S^3} X_2^\circ) \sqcup (B^4 \cup_{S^3} B^4)\] equipped with a diffeomorphism $F$ given by $f$ in the first connected component (fixing the specified $4$-ball) and the identity on the latter two components. The cyclic permutation $\tau = 312$ is even and results in \[X^\tau = (X_1 \# X_2) \sqcup (B^4 \cup_{S^3} B^4) \sqcup (B^4 \cup_{S^3} B^4)\] By Lemma~\ref{lem:connectsum}, the families Bauer-Furuta invariants of $(X,F)$ and $(X^\tau,F^\tau)$ agree. 

  Now \[\FBF^{S^1}(X^\tau,F^\tau) = \FBF^{S^1}(X_1 \# X_2, f \# \id_{X_2})\] by construction, so it remains to show that $\FBF^{S^1}(X,F) = \FBF^{S^1}(X_1,f) \wedge \BF^{S^1}(X_2)$. This will follow from careful consideration of the definition of $\FBF^{S^1}$ from the monopole map.

  In particular, $\FBF^{S^1}$ is given by pre-composition of the finite-dimensional approximation element $[\sw]$ with $p$. The key point is that, although $\FBF^{S^1}(X,\id_X) = 0$ for any $X$ (see \cite[Lemma 2.23]{lin2023isotopy}, the same proof works verbatim in the spin$^c$ setting), the monopole map for the mapping torus associated to $X$ (with fiber $X_1 \sqcup X_2 \sqcup S^4$) splits as $\SW_X = \SW_{X_1} \oplus \SW_{X_2} \oplus \SW_{S^4}$. Since the diffeomorphism on the latter two components is the identity, $\SW_{X_2}$ and $\SW_{S^4}$ can be built from the ordinary (manifold level) monopole maps by $\times \id_{S^1}$. The finite-dimensional approximation thereof is then given by \[\sw_X = \sw_{X_1} \wedge \BF^{S^1}(X_2) \wedge \BF^{S^1}(S^4)\] with the $S^1_+$ factor arising from taking the Thom space of the domain subsumed into $\sw_{X_1}$. $p^\ast$ acts by the identity on the latter two smash product factors since it only affects the $S^1_+$ factor, from which the result follows (together with the fact that $\BF^{S^1}(S^4) = [\id] \in \{*, *\}^{S^1}$, see \cite[Proposition 2.3]{bauer2004stableII}).
\end{proof}

\begin{corollary}[{cf. \cite[Remark 8.6]{krushkal2024corks}}]\label{cor:FSWfolklore}
  Exotic diffeomorphisms that are detected by the $S^1$-equivariant families Bauer-Furuta invariant cannot be isotoped to a diffeomorphism supported in a $4$-ball.
\end{corollary}
\begin{remark}
  This implies that diffeomorphisms that are detected by the families Seiberg-Witten invariant also cannot be isotopic to a diffeomorphism supported in a $4$-ball, since $\FSW$ can be recovered from $\FBF^{S^1}$ by \cite[Theorem 3.6]{baraglia2022bauer}.
\end{remark}
\begin{proof}
  Suppose $f:X \rightarrow X$ is a diffeomorphism such that $\FBF^{S^1}(X,f)$ is nontrivial, and that $f$ can be isotoped to $f'$ which is supported in a $4$-ball. Then we may write $X= X \# S^4$ where the $S^4$ summand isolates the action of $f'$, and by Corollary~\ref{cor:connectsumidentity} and isotopy invariance, \[\FBF^{S^1}(X,f) = \BF^{S^1}(X) \wedge \FBF^{S^1}(S^4,f'|_{B^4} \cup_\partial \id_{B^4})\] However, by \cite[Theorem 1.6]{lin2021family}, $\FBF^{S^1}$ vanishes for all diffeomorphisms of the $4$-sphere, so $\FBF^{S^1}(X,f) = 0$, which gives the desired contradiction.
\end{proof}

\subsection*{AI Disclosure} The first version of this paper was posted before AI tools were useful for research mathematics. An LLM was used to proofread this second version, and it suggested that Gompf's technique readily generalizes our examples to uncountable families, an extension which is included as Theorem~\ref{thm:uncountable}.



  
\setlength{\bibsep}{4pt plus 0.3ex}

{\small \footnotesize \bibliographystyle{alphamod}
\bibliography{references}}

\end{document}